\documentclass[10pt]{article}

\usepackage{graphicx}
\usepackage{amsthm}
\usepackage{amsmath,amssymb}
\usepackage{color}
\newtheorem{thm}{Theorem}
\newtheorem{cor}{Corollary}
\newtheorem{lem}{Lemma}
\theoremstyle{definition}
\newtheorem{rem}{Remark}

\newtheorem{definition}{Definition}
\numberwithin{equation}{section}
\numberwithin{definition}{section}
\numberwithin{lem}{section}
\numberwithin{thm}{section}
\numberwithin{rem}{section}
\numberwithin{figure}{section}
\title{Decomposition Approach for Low-rank Matrix Completion}
\author{Rick Ma and Samuel Cheng}

\begin{document}

\maketitle
\begin{abstract}
In this paper, we describe a low-rank matrix completion method based on matrix decomposition. 
An incomplete matrix is decomposed into 
submatrices which are filled with a proposed trimming step and then are recombined to form a low-rank completed matrix. The divide-and-conquer approach can significantly reduce computation complexity and storage requirement. Moreover, 
the proposed decomposition method can be naturally incorporated into any
existing matrix completion methods to attain further gain.
Unlike most existing approaches, the proposed method is not based on norm minimization nor SVD decomposition. This makes it possible to be applied beyond real domain and can be used in arbitrary fields including finite fields.

\end{abstract}

\section{Introduction}

Consider a large matrix with only a small portion of known entry, an interesting problem is to fill the missing entry assuming the matrix has low-rank. This problem
has several interesting applications including the so-called collaborative filtering problem \cite{srebro2004learning}. An example is the famous Netflix challenge where a huge matrix is used to represent the rating of a movie given by a user. Of course, a typical user will only rate very few movie titles. Therefore, an algorithm will be needed to complete the matrix to predict the ratings of all movies among all users. 

It has been shown theoretically that under certain assumptions the matrix can be recovered with very high accuracy \cite{candes2009exact,candes2009power,keshavan2009matrix}. Their approaches convert the rank minimization problem into a nuclear norm minimization problem instead and thus can be solved using semidefinite program (SDP). However, the complexity still grows rather rapidly with the size of the matrix $n$ ($\sim n^3$). Several efficient algorithms have been proposed including Singular Value Thresholding (SVT) \cite{SVT:08}, Atomic Decomposition for Minimum Rank Approximation (ADMiRA) \cite{ADMiRA:09}, Fixed Point Continuation with Approximate (FPCA) \cite{FPCA:09}, Accelerated Proximal Gradient (APG) \cite{APGL:09}, Subspace Evolution and Transfer (SET) \cite{dai2009set}, Singular Value Projection (SVP) \cite{SVP:09}, OptSpace \cite{keshavan2009matrix}, and LMaFit \cite{LMaFit:10}, where OptSpace and SET are based on Grassmann manifold optimization, SVT and SVP uses iterative hard thresholding (IHT) to facilitate matrix shrinkage, FPCA utilizes Bregman iterative algorithm and Monte Carlo approximate SVD, and LMaFit adopts successive over-relaxation (SOR). 

In this paper, we propose a decomposition method to allow very efficient divide-and-conquer approach when known entries are relatively very few. 
A simple ``trimming'' method is proposed to recover the decomposed ``cluster'' matrix. 
However, the decomposition method can also be combine with any other existing matrix completion techniques to yield further gain.
One advantage of the proposed approach is that unlike most existing approaches it does not utilize SVD but only relies on basic vector operations. Therefore, the approach is immediately applicable to matrices of any field (including finite field matrices). This opens up opportunities for new applications.
 

The rest of the paper is organized as follows. In the next section, we will fix our notation, describe the problem precisely, and present several properties to be used in the later sections. Sections \ref{sect:diag} and \ref{sect:sdiag} will describe the decomposition procedures and present our main results. Section \ref{sect:trimming} will describe the trimming process.

%
%

\section{Minimum Rank of Incomplete Matrix}

Let us start with a few notes on our notation.
When things are clear, lines of partition in matrices will not be shown; the $?$ sign may represent an unknown entry, a row or column of unknown entry, a matrix of unknown entry, etc; and similar for the $0$ sign.

Given a finite size matrix $M$ over field $\mathbb{F}$ to be completed, let
\begin{equation}
 S(M)=\{\bar M | \bar M \mbox{ is a completion of } M\}. 
\label{rick:1}
\end{equation}

If $ M$ is already completed, then $ S(M)=\{M\}$.
We define
\begin{equation}
mr(M)=\min_{\bar M \in S(M)} \mbox{rank} \bar M. 
\label{rick:2}
\end{equation}

Such minimum exists 
 because $ \mbox{rank} (S(M)) \subset \mathbb{N}$ and hence $ \exists \bar M \in S(M)$ such that 
\begin{equation}
\mbox{rank} \bar M= mr M.
\label{rick:new1.3}
\end{equation}
If $M= \begin{pmatrix} A & B \\ C & D \end{pmatrix}$,
\begin{equation}
 \exists \bar A \in S(A) \mbox{ such that } mr(M) = mr 
\begin{pmatrix} \bar A & B \\ C & D \end{pmatrix}, 
\label{rick:new1.4}
\end{equation}
as we can always find $\bar A$ from $\bar M$ in \eqref{rick:new1.3}. We list other properties about $mr(M)$ that will be quoted:
\begin{align}
 &mr (M) \le \mbox{rank} \bar M,  \qquad \forall \bar M\in S(M)       
\label{rick:3} \\ 
& mr(M) \le mr(P) \mbox{ if  $P$ is any partial completion of  $M$ } 
\label{rick:4} \\ 
 &mr([A | B]) \le mr A + mr B    
\label{rick:5}  \\ 
& mr M^t = mr M  
\label{rick:6} \\ 
 &mr \left( \left[\begin{matrix} A & B \\ C & D \end{matrix}\right]\right) \ge mr(A)
\label{rick:7} \\
&mr M = mr N \mbox{ if $ N$ can be obtained from $ M$ through
interchanging of columns/rows.} 
\label{rick:8} 
\end{align}


\subsection{Junk Row and Junk Column}
\begin{definition}
A row(column) contains entirely either zero or unknown  will be refered as a junk row(column).
\label{def:9}
\end{definition}
Certainly, we have 
\begin{equation}
\label{rick:10}
 mr(J)=0 \mbox { if }  J \mbox{ is a junk row (column), }                   
 \end{equation}
since we can always complete $J$ entirely by zero entries.


\begin{thm}
\label{thm:1}
Let $ M=[J| N]$ where $ J$ is a junk column, then
$ mr M= mr N$.
\end{thm}

\begin{proof}
By \eqref{rick:5}, 
we have
$ mr M \le mr J + mr N = mr N$.
On the other hand $ mr M \ge mr N$ by \eqref{rick:7}. 
Hence
$ mr M=mr N$.
\end{proof}
%

Thanks to \eqref{rick:6}, we have the following corollary:
\begin{cor}
\label{cor:1}
$ mr \left( \left[ \begin{matrix} J \\ N \end{matrix} \right] \right) = mr N$ if $ J$ is a junk row.
\end{cor}


\subsection{Equivalence}
We say $ M$ is equivalent to $ N$ and write $ M \sim N$ iff $ N$ can be
obtained from $ M$ through row interchanging, column interchanging, 
and junk rows and columns deletion and augmentation.
By Theorem \ref{thm:1} and 
and \eqref{rick:8}, we have
\begin{equation}
 M \sim N \Rightarrow mr M = mr N.                                    
\label{rick:11}
\end{equation}

\section{Unknown-diagonalization}
\label{sect:diag}
Define 
\begin{equation}
\mbox{u-diag}(B_1,B_2,\cdots,B_n)\triangleq \left[\begin{matrix} 
B_1 & \mbox{\bf ?} &\cdots & \mbox{\bf ?}\\
\mbox{\bf ?} &B_2&\cdots&\mbox{\bf ?}\\
\cdots&&\cdots\\
\mbox{\bf ?}&\mbox{\bf ?}&\cdots&B_n\end{matrix}\right]. \label{rick:12}
\end{equation}
We say $M$ is $\mbox{u-diagnonalizable iff } M \sim \mbox{u-diag} (A,B)$,
and both $ A$ and $B$ contain at least one nonzero known entry.

\begin{thm}
Let $ M\sim \mbox{u-diag}(B_1,\cdots,B_n)$, then
$ mr M= max_{1\le i\le n} mr (B_i)$. 
\label{thm:2}
\end{thm}

\begin{proof}
By \eqref{rick:11} and induction, all we need to show is when
$ M=\mbox{u-diag}(A,B)$.
Let $ mr A=a$ and $ mr B=b$, by \eqref{rick:7} we have
\begin{equation}
 mr M \ge max (a, b) \label{rick:13} 
\end{equation}

Let $ \bar A$, $ \bar B$ be completions of $ A$ and $ B$, respectively
such that
$ \mbox{rank} \bar A= a$ and $ \mbox{rank} \bar B=b$ (c.f. \eqref{rick:new1.3}), then by \eqref{rick:4},
\begin{equation}
mr M\le mr \left( \mbox{u-diag}(\bar A, \bar B) \right) 
\label{rick:14} 
\end{equation}

Combining \eqref{rick:13} and \eqref{rick:14}, we conclude that
\begin{equation}
 max(a,b) \le mr M \le mr \left( \mbox{u-diag}(\bar A, \bar B) \right).
\label{rick:15} 
\end{equation}

By \eqref{rick:8}, we can simply assume the first $ a$ columns
of $ \bar A$ form a basis of $ Col \bar A$ ; and the first $ b$ columns
of $ \bar B$ form a basis of $ Col \bar B$.
Without loss of generality, let us assume $ a\ge b$. We complete the
﻿﻿matrix $ \mbox{u-diag}(\bar A,\bar B)$ by filling up the columns: 
\begin{align}
\left[ \begin{matrix}
        \bar A_i \\ \mbox{\bf ?}
       \end{matrix}
\right] & \Rightarrow
\left[ \begin{matrix}
        \bar A_i\\ \bar B_i
       \end{matrix}
\right], 
\nonumber
\\
\left[ \begin{matrix}
        \mbox{\bf ?} \\ \bar B_i 
       \end{matrix}
\right] & \Rightarrow
\left[ \begin{matrix}
        \bar A_i\\ \bar B_i
       \end{matrix}
\right], & \mbox{for $1 \le i \le b$},
\label{rick:new2.5}
 \\
\left[ \begin{matrix}
        \bar A_i \\ \mbox{\bf ?}
       \end{matrix}
\right] & \Rightarrow
\left[ \begin{matrix}
        \bar A_i\\ \bar 0
       \end{matrix}
\right], & \mbox{for $b+1 \le i \le a$}.
\label{rick:new2.6}
\end{align}
For $i > a$, we make use of the fact that $\bar A_i$ is a linear combination of 
$\{ \bar A_k | k \le a \}$. We fill
\begin{align}
\left[ \begin{matrix}
        \bar A_i \\ \mbox{\bf ?}
       \end{matrix}
\right] & \Rightarrow
\left[ \begin{matrix}
        \bar A_i\\ \sum_{k=1}^a a_{i,k}\bar B_k
       \end{matrix}
\right], & \mbox{for $i > a$}
\label{rick:new2.7}
\end{align}
Similarly,
\begin{align}
\left[ \begin{matrix}
        \mbox{\bf ?} \\ \bar B_i
       \end{matrix}
\right] & \Rightarrow
\left[ \begin{matrix}
       \sum_{k=1}^b b_{i,k}\bar A_k \\   \bar B_i 
       \end{matrix}
\right], & \mbox{for $i > b$}, 
\label{rick:new2.8} 
\end{align}



Now we have a completed $ \mbox{u-diag}(\bar A,\bar B)$ and the first $ a$ of its columns
$ \left[\begin{matrix}\bar A_1\\\bar B_1\end{matrix}\right]$, $ \left[\begin{matrix}\bar A_2\\\bar B_2\end{matrix}\right]$,$ \cdots$,$ \left[\begin{matrix}\bar A_b\\\bar B_b\end{matrix}\right]$,$ \left[\begin{matrix}\bar A_{b+1}\\{\bf 0}\end{matrix}\right]$,$ \cdots$,$ \left[\begin{matrix}\bar A_a\\{\bf 0}\end{matrix}\right]$ form a basis for its column space. Hence it has rank $a$. By \eqref{rick:3}, $ mr \left( \mbox{u-diag}(\bar A,\bar B) \right)\le a$ and hence by \eqref{rick:15} and $ max(a,b)=a$, we get $ mr M=a=mr \left( \mbox{u-diag} (A,B) \right)$ as wanted.
\end{proof}

\begin{rem}
\label{rem:2.1}
Suppose $A$ has been completed by $\bar A$ and $a= \mbox{rank}(\bar A)$. Then if the number of column of $B=n \le a$, then we can complete $B$ arbitrarily and do the completion in \eqref{rick:new2.5}-\eqref{rick:new2.8} as if $b=n$. More generally, given $M=\mbox{u-diag}(\bar A,B)$ with $\bar A$ is completed. Then the completing process of $B$ can be stopped once we know that the final $\mbox{rank} B$ will not be greater than $\mbox{rank} \bar A$ no matter how we do the remaining completion on $B$. For example, if $\mbox{size}(B) \le \mbox{rank} \bar A$, where
\begin{align}
 \mbox{size}(B) = \min (\mbox{number of column of $B$}, \mbox{number of row of $B$}), 
\label{rick:new2.9}
\end{align}
%
then we can complete $B$ arbitrarily to start with.
\end{rem}

\subsection{Percolation and Clusters}

\begin{figure}
\begin{center}
 \includegraphics[width=5in]{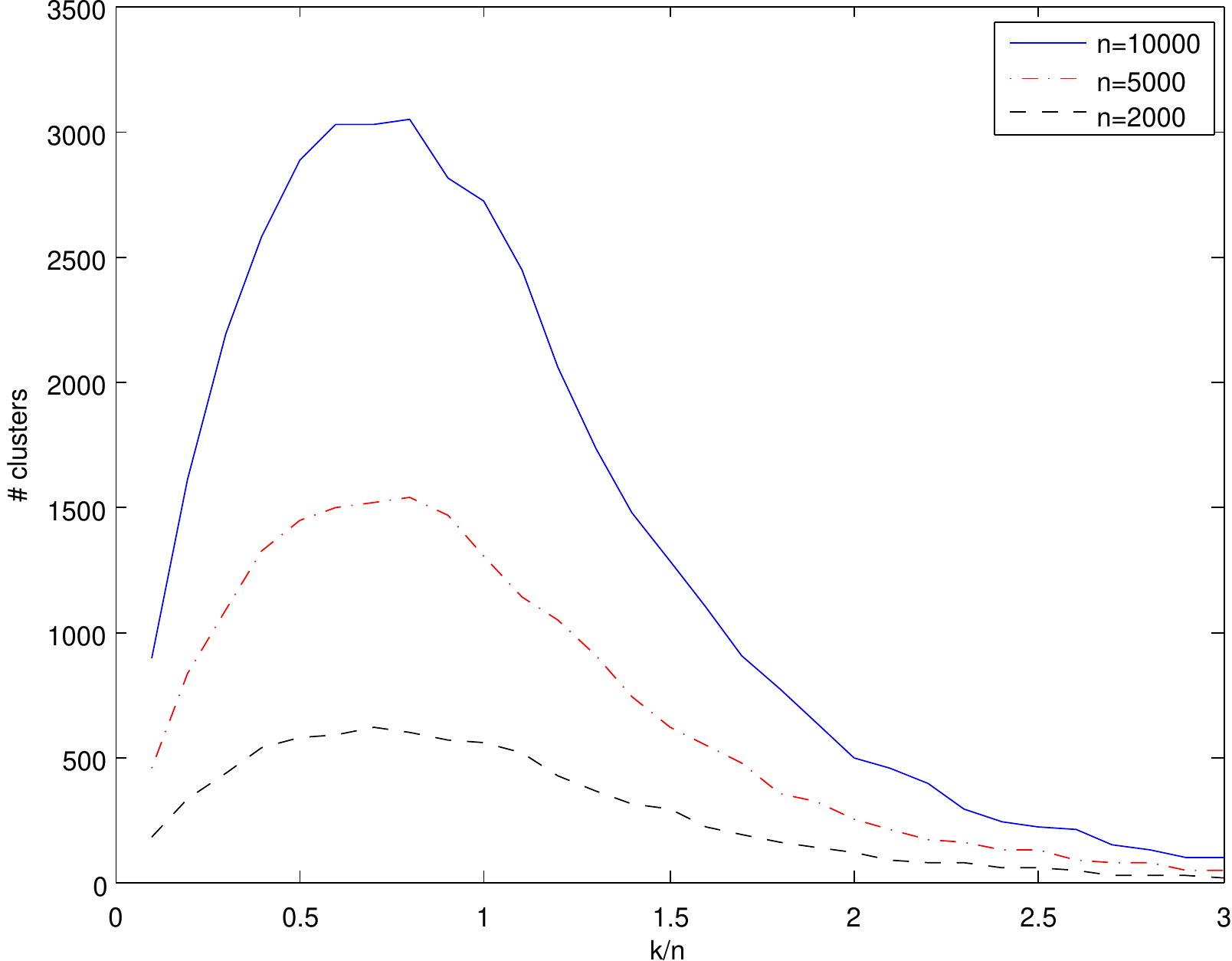}
\caption{The number of clusters versus $k/n$, where $k$ is the number of known entries and the size of the matrix is $n \times n$.}
\label{fig:clusters}
\end{center}
\end{figure}

We would call those $B_1,\cdots B_n$ in Theorem \ref{thm:1} as clusters. In other words, clusters are matrices that cannot be u-diagonized. They are not the clusters in the 2-d square lattice, where each point, not counting the edgy one, has 4 neighbors. In our case, each entry in an $n \times m$ matrix has $n + m - 2$ neighbors, from the view of percolation. Despite that difference, the two models share the same percolation threshold at $p \sim 0.6$ \cite{stauffer1994introduction}, where $p$ is the occupation rate. That means if about $60 \%$ of our entries are known, then there is probably one cluster left and the matrix cannot be u-diagonalized. We estimate the number of clusters as the size of the matrix and the number of known entries vary through Monte Carlo simulation and the results are shown in Figure \ref{fig:clusters}. We can see the number of clusters increases as as the number of known entries increases and is peak when the occupation ratio is at about $0.7$ regardness the size of the matrix.

\subsection{A Decomposition Algorithm} 
\label{sect:2.2}

First of all, we set all junk rows and junk columns of the given matrix to zero and blackout them. Now we are working with a junk-free matrix.

We create a row set and a column set for the first cluster. Then put the row position of the first row into the row set; the columns positions of the columns with known entries in that row into the column set. Thanks to \eqref{rick:8}, sorting is not necessary. We black out the row against repeating searches. For each new comers of the column set, we search vertically for its known entries and put the corresponding row position into the row sets. After that we black out the searched columns. Now the row set may have new comers. We enlarge the column set in the way that we enlarged the row set. Both sets keep growing until there is no more new comer.

Then we create another column set and another row set. Repeat the procedure for the next cluster, if the remaining matrix has not been blackout to void.

\section{Sub Unknown-diagonalization}
\label{sect:sdiag}
\begin{definition}
\label{definition:3.1}
Let $u$ be a vector and $A$ be a matrix. We define 
 $zero(u, A)$ as $1$ if  $\forall [\bar u | \bar A]\in S([u | A])$
 (c.f. \eqref{rick:1})
with $\mbox{rank} [\bar u | \bar A]= mr [u | A]$, then $\bar u= {\bf 0}$.
Otherwise, we define $zero(u, A)$ as $0$.
\end{definition}


\noindent
E.g.
$
zero \left( 
\begin{pmatrix} 0 \\ {\bf ?} \end{pmatrix},
\begin{pmatrix} 1 \\ 1 \end{pmatrix}
\right)=1$,
$zero\left(
\begin{pmatrix}0 \\ {\bf ?} \end{pmatrix},
\begin{pmatrix}0 \\ 1 \end{pmatrix}
\right)=0$.
{\color{black} Notice that the necessary condition (not sufficient) for $zero (u, A)=1$ is
that $u$ must be a junk column.}
Hence we have 
$zero(u,A) \le junk(u)$, where $junk(u)=\left\{\begin{array}{l l}
1, &               \mbox{if $u$ is junk,} \\
0, & otherwise.                           
                                        \end{array}\right.$.


\begin{thm}
\label{thm:5}
Let $M= 
\begin{pmatrix}
v_1 & B1& &        \mbox{\bf ?} \\
v_2  & &     B2  \\
\cdots & & &           \cdots \\
v_n &     \mbox{\bf ?} &&          Bn 
\end{pmatrix}
$, where $B_i$ are matrices, $v_i$ are vectors that
$\begin{pmatrix}
v_1 \\
v_2 \\
\cdots \\
v_n  
 \end{pmatrix}
$ is not a junk column, then
$mr M= \max_i (mr[v_i | B_i] + zero (v_i, Bi))$.
\end{thm}

\begin{proof}
We will show the case when $n=2$; cases of higher $n$ are easy to generalized. Let
$M=\begin{pmatrix} u & A & \mbox{\bf ?} \\
v &\mbox{\bf ?} & B\end{pmatrix}$ with 
$\begin{pmatrix}
  u \\v 
 \end{pmatrix}
$ is not a junk column.
If $zero (u, A)=0$, then we have
\begin{align}
mr M & \ge mr [u | A]  \nonumber \\
& =mr [u | A] + zero (u, A), \label{rick:25}
\end{align}
where the first inequality is by \eqref{rick:7}.
If $zero (u, A)=1$, then $u$ must be a junk column. Since
$\begin{pmatrix}
u \\
v  
 \end{pmatrix}
$ is not a junk column, $v$ must contains a nonzero known
entry $\lambda$.
Let $M' = \begin{pmatrix}
u   &         A\\
\lambda  & {\bf ?} \cdots {\bf ?}
          \end{pmatrix}
$. By \eqref{rick:5} and \eqref{rick:7},
\begin{equation}
 1+mr A \ge mr M' \ge mr A. \label{rick:26}
\end{equation}
So there are only two possibilies for $mr M'$.
Assume $mr M'= mr A$. By 
\eqref{rick:new1.3}
we can
pick $[\bar u | \bar A]\in S([u | A])$ such that
\begin{equation}
mr A= mr M' = mr \begin{pmatrix} \bar u   &   \bar A \\
\lambda   & {\bf ?}\cdots {\bf ?}\end{pmatrix}
\ge 
mr [ \bar u | \bar A ].
\label{rick:27} 
\end{equation}
Therefore, 
$mr A =mr [\bar u \bar A]$ as
$rank[\bar u \bar A]\ge rank \bar A \ge mr A$. Since $zero (u,A)=1$,
we have
$\bar u={\bf 0}$ (c.f. def \ref{definition:3.1}).
But that will make
$mr M’ = rank \bar A +1$ because $\lambda \neq 0$.
We conclude that $mr M'$ cannot equal to $mr A$ and hence we must have
\begin{align}
mr M & \ge mr M’= mr A+1  &\mbox{(c.f. \eqref{rick:26})} \nonumber \\
     & =  mr [u | A] + zero (u, A). & \mbox{(c.f. Theorem \ref{thm:1})}   
\label{rick:new3.4}      
\end{align}

Combine \eqref{rick:25} and \eqref{rick:new3.4} and the symmetry between $(u, A)$
and $(v, B)$ we get
\begin{equation}
mr M \ge \max (mr[u | A] +zero (u, A), mr[v | B] + zero (v, B)). 
\label{rick:new3.5} 
\end{equation}

Without loss of generality, we may assume
\begin{align}
a=mr[u| A] + zero (u,A) \ge mr [v| B] + zero (v, B)=b.           
\label{rick:new3.6} 
\end{align}

Pick $[\bar A_0| \bar A]\in S([u| A]), [\bar B_0| \bar B]\in S([v| B])$
s.t. 
\begin{align}
rank  [\bar A_0| \bar A]& = mr [u | A], \nonumber \\
 rank  [\bar B_0| \bar B]& = mr [v| B].  \qquad \mbox{(c.f. \eqref{rick:2})} 
\label{rick:new3.7}
\end{align}
%
%
and 
\begin{align}
\bar A_0 &\neq  \vec 0,  & \mbox{ if $zero (u, A)=0$},& \nonumber \\
  \bar B_0 &\neq \vec 0, & \mbox{ if $zero (v, B)=0$}. 
&\qquad  \mbox{(c.f. def \ref{definition:3.1})}     
\label{rick:new3.8}
 \end{align}

Thanks to \eqref{rick:8}, we can assume
$$Span \{\bar A_0, \bar A_1,....,\bar A_{a-1}\}=Col [\bar A_0| \bar A]$$,
$$Span \{\bar B_0, \bar B_1,....,\bar B_{b-1}\}=Col [\bar B_0| \bar B]$$.

Then we complete the rest by filling up the columns:
\begin{align}
\left[ \begin{matrix}
        \bar A_i \\ \mbox{\bf ?}
       \end{matrix}
\right] & \Rightarrow
\left[ \begin{matrix}
        \bar A_i\\ \bar B_i
       \end{matrix}
\right], 
\nonumber
\\
\left[ \begin{matrix}
        \mbox{\bf ?} \\ \bar B_i 
       \end{matrix}
\right] & \Rightarrow
\left[ \begin{matrix}
        \bar A_i\\ \bar B_i
       \end{matrix}
\right], & \mbox{for $1 \le i < b$},
\label{rick:new3.9}
 \\
\left[ \begin{matrix}
        \bar A_i \\ \mbox{\bf ?}
       \end{matrix}
\right] & \Rightarrow
\left[ \begin{matrix}
        \bar A_i\\ \bar 0
       \end{matrix}
\right], & \mbox{for $b \le i < a$}.
\label{rick:new3.10}
\end{align}
For $i \ge a$, we can pick $\{a_{i,k} \} \in \mathbb{F}$ s.t. 
$\bar A_i = \sum_{0 \le k < a} a_{i,k} \bar A_k$ and fill the columns
\begin{align}
\left[ \begin{matrix}
        \bar A_i \\ \mbox{\bf ?}
       \end{matrix}
\right] & \Rightarrow
\left[ \begin{matrix}
        \bar A_i\\ \sum_{0 \le k < a} a_{i,k}\bar B_k
       \end{matrix}
\right], & 
\label{rick:new3.11}
\end{align}
Similarly, for $i \ge b$, we pick $\{ b_{i,k} \} \in \mathbb{F}$ s.t. 
$\bar B_i = \sum_{0 \le k < b} b_{i,k} \bar B_k$ and fill the columns
\begin{align}
\left[ \begin{matrix}
        \mbox{\bf ?} \\ \bar B_i
       \end{matrix}
\right] & \Rightarrow
\left[ \begin{matrix}
       \sum_{0 \le k < b} b_{i,k}\bar A_k \\   \bar B_i 
       \end{matrix}
\right], & 
\label{rick:new3.12} 
\end{align}
Then the completed matrix is rank $a$ with the first $a$ columns form a basis of its column space. By \eqref{rick:3}, 
\begin{align*}
      mr M & \le  a  &\\
          &= rank [\bar A_0 \bar A] + zero (u, A)  & 
\mbox{(c.f. \eqref{rick:new3.6}, \eqref{rick:new3.7})} \\
           &=\max (mr[u A] +zero (u, A), mr[v B] + zero (v, B)) &
\end{align*}
as assumed. Together with \eqref{rick:new3.5}, we get
$mr M=\max(mr[u A] +zero (u, A), mr[v B] + zero (v, B))$.

\end{proof} 

\begin{rem}
\label{rem:3.1}
 Suppose $[u| A]$ has been completed by $[\bar u| \bar A]$ and we have $a=zero (\bar u, \bar A) + rank[\bar u \bar A]$. Then if the number of column of $[v | B]=n \le a$, we can complete $[v| B]$ arbitrarily and do the completion in \eqref{rick:new3.9}-\eqref{rick:new3.12}
 as if $b=n$. More generally, the completing process of $B$ can be stopped once we know that the final $zero(\bar v, \bar B) + rank [\bar v \bar B]$ won't be greater than $a$ no matter how we do the remaining completion on $[v B]$. For example if $junk(v)+size([v B]) \le a$ (c.f. \eqref{rick:new2.9}, def \ref{definition:3.1}), then we can complete $[v| B]$ arbitarily at the beginning and do the completion \eqref{rick:new3.9}-\eqref{rick:new3.12}.
\end{rem}


\subsection{How to decompose sub u-diagonalizable matrix}

\begin{definition}
 A matrix $C$, not u-diagonalizable, becomes u-diagonalizable after deleting a row or a column is called sub u-diagonalizable. The row (column) is called conjoined row (column).  
\end{definition}

For example, 
$\begin{pmatrix}
  v_1 \\ v_2 \\ \cdots \\v_n
 \end{pmatrix}
$
in Theorem \ref{thm:5} is a conjoined column and $M$ becomes 
u-diagonalizable without it.

 \begin{definition}
\label{definition:3.3}
Given two vectors $v$ and $w$ of same length, we say $v$ is a {\em donor} for $w$ $(v \succeq w)$ iff all of the unknown positions of $v$ are also unknown in $w$. In order words, after some row interchanging $[w | v]
= 
\begin{pmatrix}
 \bar r & \bar d \\
\vec ?  &  n
\end{pmatrix}
$ with $\bar r$ and $\bar d$ are completed. Clearly, $v \succeq  w$ and $w \succeq u$ imply $v \succeq u$. However, if $v \succeq w$ and $w \succeq v$, we do not have $v=w$. Vectors $v$ and $u$ are said to be {\em comparable} if either $v \succeq w$ or $w \succeq v$.    
 \end{definition}

\begin{thm}
\label{thm:3.2}
 Conjoined row (column) does not have donors among other rows of the sub u-diag matrix.
\end{thm}

\begin{proof}
 Let $C$ be the sub u-diag  matrix, then it must have the following structure (after some row and column interchanging):  
$
C = \begin{pmatrix}
 A & ? \\
           ? & B  \\
           u & v
\end{pmatrix}, 
$
where $u$ cannot be entirely unknown, otherwise $C$ is u-diagonalizable. Now, rows in $[? | B]$ cannot be donors of $[u | v]$, the conjoined row. Similarly $v$ cannot be entirely unknown and hence, rows in $[A | ?]$ cannot be donors of $[u| v]$ neither.  
\end{proof}

Therefore if $C$ is a sub unknown-diagonalizable, we will not miss the chance of decomposing it if we have tested every row and column that does not have a donor. That is to blackout the suspicious row (column) and then carrying out the decomposition mentioned in Section \ref{sect:2.2}. We would like to call the decomposed components as {\em sub-clusters}. For example,
$\begin{pmatrix}
 A \\u
\end{pmatrix}$
and $\begin{pmatrix}
     B \\v
    \end{pmatrix}$
are sub-clusters of the $C$ in the proposition.  

Unlike cluster that cannot be further unknown-diagonalized, sub-clusters can be sub unknown-diagonalizable. For example
$C=\begin{pmatrix}
A & ? & ?  \\     
 u & v & ?   \\ 
? & B & ?    \\
? & x & y     \\
? & ? & D  
 \end{pmatrix},$
where both $[u | v | ?]$ and $[?| x| y]$ are conjoined rows. In that case, we may first decompose $C$ into sub-clusters
$\begin{pmatrix}
  A \\u
 \end{pmatrix}$
and 
$\begin{pmatrix}
  v & ? \\ B & ? \\x & y \\ ? & D
 \end{pmatrix}$.
 Then we may further decompose the later into
$\begin{pmatrix}
  v \\ B \\ x 
 \end{pmatrix}$
and 
$\begin{pmatrix}
  y \\ D
 \end{pmatrix}$, if necessary.

\section{Trimming}

\label{sect:trimming}
\begin{lem}
 $mr [v | M]=mr M$ if $\forall \bar M \in S(M)$,
$Col (\bar M) \cap S(v) \neq \mathbb{\phi}$,.
\label{lem:4.1}
\end{lem}
\begin{proof}
 By (1.9) we already have
$mr[v| M] \ge mr M$.                   
Let $\bar M \in S(M)$ such that $rank \bar M=mr M$ (c.f. \eqref{rick:2}).
Then pick a $\bar v \in S(v) \cap Col(\bar M)$. From \eqref{rick:3}, we get
$mr [v | M] \le rank [\bar v | \bar M]=rank \bar M=mr M$.       
\end{proof}

\begin{thm}
\label{thm:4}
Let $M_{d_i}$ be the $d_i$-th column of $M$ and donor (c.f. def \ref{definition:3.3})
of a vector $v$ for $1\le i \le t$, such that after some row interchanging,
$[v | M_{d_1} M_{d_2} \cdots M_{d_t}]=
\begin{pmatrix}
\bar r | \bar D \\
                             ?  | N    
\end{pmatrix}$
 with $\bar r$ and $\bar D$ are completed.
If $\bar r\in Col (\bar D)$, then $mr[v | M]= mr M$. 
\end{thm}

\begin{proof}
Thanks to \eqref{rick:8}, we may start with
$
[v |M]=
\begin{pmatrix}
 \bar r |  \bar D A \\
? | N B
\end{pmatrix}
$
%

Pick $a_i\in \mathbb F$ such that 
\begin{equation}
\bar r=\sum a_i \bar D_i. 
\label{rick:new4.1} 
\end{equation}
Then $\forall
\bar M=
\begin{pmatrix}
 \bar D \bar A \\
\bar N \bar B 
\end{pmatrix}\in S(M)$,
we complete $v$ to 
\begin{align}
 \bar v=
\begin{pmatrix}
\bar r \\
        \sum a_i\bar N_i 
\end{pmatrix}
\in S(v)\cap Col (\bar M). 
\label{rick:new4.2}
\end{align}
(Note that $rank[\bar v | \bar M]= rank \bar M$.)
Now Lemma \ref{lem:4.1} implies $mr[v | M]=mr M$.
 \end{proof}

\subsection{Trimming Process}

We test column by column to see if we can make use of Theorem \ref{thm:4} to trim away some columns from a given matrix, which is probably a sub-cluster mentioned in the previous section. We call this process as {\em column trimming}.  
When we find a column satisfying the condition of Theorem \ref{thm:4}, we will mark down the dependency relation between it and its donor (i.e. 
\eqref{rick:new4.1}) in order. Then we black it out and go for the next column. 

Similarly, we have row trimming. An uninterrupted (c.f. Remarks 
\ref{rem:2.1} and \ref{rem:3.1}) trimming process starts with a column trimming followed by a row trimming, or the other way round. Then we carry out these two kinds of trimming one after the other, until there is no more reduction in the matrix. 
After the trimmed matrix gets completed, we restore, in reverse order, the blackouts with the completed forms given by \eqref{rick:new4.2}.  

The following proposition is interesting in its own right and may be useful for our future study.

\begin{lem}
\label{lem:4.2}
 If $v$ is a column that $S(v)\cap Col(\bar M)=\mathbb{\phi}$, $\forall
\bar M\in S(M)$. Then $mr[v  | M]=mr M+1$.
\end{lem}
 
\begin{proof}
By \eqref{rick:5} $mr [v| M] \le mr M +1$. Let us complete $[v | M]$ by $[\bar v | \bar M]$ such that
$rank [\bar v | \bar M]= mr [v | M]$ (c.f. \eqref{rick:2}). Since $\bar v \notin Col (\bar M)$, we must have $rank [\bar v | \bar M]= rank \bar M+1$, which implies
$mr [v | M] = rank \bar M + 1 \ge mr M + 1$ (by \eqref{rick:3}).  
\end{proof}


\noindent
{\bf Proposition.} {\it
Suppose $M$ is a matrix that every two columns $M_i$ and $M_j$ of $M$ are comparable (c.f. def \ref{definition:3.3}). Then one round of column trimming followed by arbitrarily completion and proper restoration (i.e. \eqref{rick:new4.2}) of the trimmed columns complete $M$ to its minimum rank.  }

\begin{proof}
Let $T$ be the trimmed $M$. Every two columns of $T$ are also columns of $M$ and hence comparable. So there exists a column $T_j$ in $T$ such that $T_j \preceq T_i$ for all $i$. After some row interchanging and column interchanging, we have
$T_j=\begin{pmatrix}
  \bar r \\? 
 \end{pmatrix}$
and 
$T = \begin{pmatrix}
      \bar r & \bar D  \\
? & N
     \end{pmatrix}$
with $\bar r$ and $\bar D$ are completed. Notice that $\bar r\notin Col (\bar D)$, otherwise $T_j$ has been trimmed away already. Now Lemma \ref{lem:4.2} and \eqref{rick:5} give $mr N=1+ mr
\begin{pmatrix}
 \bar D \\ N
\end{pmatrix}
$. Repeating the argument, we get
$mr N= \mbox{number of columns of }N \ge rank (\bar N)$, $\forall \bar N\in S(N)$. Together with Theorem 
\ref{thm:4} and 
\eqref{rick:3}, we conclude that $mr M=mr N= rank (\bar N)$, $ \forall N\in S(N)$. Finally, the proper restoration (c.f. \eqref{rick:new4.2}) restores $\bar N$ to $\bar M\in S(M)$ such that $rank (\bar M)=rank (\bar N)=mr (M)$.
 \end{proof}

\subsection{Trimming Process with Approximation}
The trimming process stops when there are no more columns or rows fulfill the condition of Theorem \ref{thm:4}. But we can always make an approximation by blacking out a column or a row as if it fulfill the condition and continue the trimming. With this approximation, the process stops when there is no more unknown left in the trimming matrix. (We may even choose the row or column that has no donor (c.f. Thm \ref{thm:3.2}) to black out and check for u-diagonalization. )  

Then we restore the blackouts in reverse order. When we meet a blackout without dependency relation (i.e. \eqref{rick:new4.1}) to restore, we check for the condition of Theorem \ref{thm:4}, again. The first time was with the uncompleted trimming matrix; this time is with the completed restoring matrix. If the condition is fulfilled, we restore the blackout with completed form given by
\eqref{rick:new4.2}. This will not compromise (further) the minimum rank that we can reach. Otherwise we restore the blackout with arbitrary completed form, which may cause one (more) rank deviation from the possible minimum.

\bibliographystyle{IEEEtran}
\bibliography{../../ref}
\end{document}